\documentclass[10pt,notitlepage,oneside,openany]{amsart}
\usepackage{amssymb}
\usepackage{amsthm}
\usepackage{booktabs}
\usepackage{mathpple}

\theoremstyle{plain}

\newtheorem*{theo}{Theorem}

\numberwithin{equation}{section}

\def\eq{\quad\Longrightarrow\quad}

\begin{document}

\begin{center}
\textsc{\Large\bf Divide-and-conquer generating functions}\\[5pt]
{\Small Part I}\\[5pt]
{\bf Elementary Sequences}\\[10pt]
Ralf Stephan\footnote{\texttt{mailto:ralf@ark.in-berlin.de}}\\[15pt]
\begin{minipage}[t]{5in}\parindent4mm \baselineskip8.5pt \Small
Divide-and-conquer functions satisfy equations in $F(z),F(z^2),F(z^4)\ldots$.
Their generated sequences are mainly used in computer science, and they
were analyzed pragmatically, that is, now and then a sequence was picked
out for scrutiny. By giving several classes of ordinary generating functions 
together with recurrences, we hope to help with the analysis of many such 
sequences, 
and try to classify a part of the divide-and-conquer sequence zoo.
\end{minipage}
\end{center}

\bigskip

Nowadays, it is routine to work within the equivalence of linear recurrences
and rational o.g.f.s., or even broader, with hypergeometrical functions and
binomial identities. Mysteries, however, remain with many
sequences generated by recurrences of other type. Empirical studies played
the main part of this work where we focus on recurrences of the form
\begin{equation}\label{rec}
a_{0\text{\ or\ }1}=c, \quad 
a_{2n}=f(a_n,a_{n+i_1},\ldots,n), \quad 
a_{2n+1}=g(a_n,a_{n+j_1},\ldots,n),
\end{equation}
and their power series generating functions.

In the first section, we introduce the reader to a class of functions that
can generate said sequences. Then, several g.f.s are given together with
corresponding recurrences and, partly, informal proofs. The last two sections
list example sequences and some further questions.

\section{Introduction and Definitions}

The first wide survey of the sequences generated by recurrences of form \eqref{rec} was
done by Allouche and Shallit\cite{b-kreg} who managed to give a broad definition
of such sequences by the notion of 2-regularity. As many of the k-regular
sequences are interesting to computer scientists, but k-regularity seems too
fuzzy a concept, it was necessary to find more detailed properties that
can be used to both distinguish between sequences and make them more amenable
to analysis. One fruitful attempt in this regard is the fitting term of
\emph{divide-and-conquer} by Dumas who revamped\cite{b-dumas93} an old concept by 
Mahler---Dumas contemplates\cite{b-dumas-www} the following definitions.

A formal series $f(z)$ is called {\bf Mahlerian} if it satisfies an
equation
\begin{equation}
c_0(z)f(z)+c_1(z)f(z^2)+\dotsb+c_N(z)f(z^{2^N})=0
\end{equation}
where the $c_k(z)$ are polynomials, not all zero. Likewise,
Mahlerian sequences satisfy a recurrence of the form
\begin{equation}
\sum_{\ell}c_{0,\ell}f_{n-\ell}+\sum_{\ell}c_{1,\ell}f_{(n-\ell)/2}+\dotsb
+\sum_{\ell}c_{N,\ell}f_{(n-\ell)/2^{\ell}}=0,
\end{equation}
where the sums are finite and the coefficients $c_{k,\ell}$ are not
all zero.

A series $f(z)$ is of the {\bf divide-and-conquer }
(DC) type  if it satisfies an equation with a right member
\begin{equation}
c_0(z)f(z)+c_1(z)f(z^2)+\dotsb+c_N(z)f(z^{2^N})=b(z)
\end{equation}
in which $b(z)$ is a formal series and $c_k(z)$
are polynomials not all zero.  A sequence is of the divide-and-conquer
type if its generating series is of the divide-and-conquer type.

\bigskip

{\bf Examples.} The sequence $\nu(n)=e_1(n)$, which counts ones in the
binary representation of~$n$, is of DC type since its g.f.~is
\begin{align*}
F(z) &=\sum_{i\ge0}e_1(i)z^i=\frac1{1-z}\sum_{k\ge0}\frac{z^{2^k}}{1+z^{2^k}}
\quad = \quad \frac1{1-z}\left(\frac{z}{1+z}+\frac{z^2}{1+z^2}+\cdots\right),\\
F(z^2) &= \frac1{1-z^2}\left(\frac{z^2}{1+z^2}+\frac{z^4}{1+z^4}+\cdots\right) 
=\frac1{1-z^2}\left((1-z)F(z)-\frac{z}{1+z}\right),\\
&\Longrightarrow \qquad (1-z^2)F(z^2)-(1-z)F(z)=-\frac{z}{1+z}.
\end{align*}
The Thue-Morse sequence on $\{1,-1\}$, $t(n)$, is not only of DC type but also Mahlerian:
\begin{align*}
T(z) &=\sum_{i\ge0}t(i)z^i=\prod_{k\ge0}1-z^{2^k}=(1-z)(1-z^2)(1-z^4)\cdots,\\
T(z^2) &=(1-z^2)(1-z^4)\cdots\quad=\quad T(z)/(1-z),\\
&\Longrightarrow \qquad T(z^2)-T(z)/(1-z)=0.\qquad\qquad\qquad\qquad\qquad\qquad\qquad\qquad\blacklozenge
\end{align*}

We will call generating functions and the associated sequences {\bf of
elementary divide-and-conquer type}, if 
\begin{itemize}
\item the g.f.~can be written as
sum or product of rational functions and  
infinite sums or products in the unknown (here~$z$), 
\item where~$z$ appears in the sum/product only as a rational
function of~$z^{2^k}$ (where~$k$ is the sum/product index starting with~$0$). 
\item The index~$k$ may appear in the product/sum  additionally as
exponent of an integer factor, and we allow only integer coefficients
to the monomials. 
\end{itemize}
The previous examples are both elementary in this
regard.

That the above defined functions are of DC type is easily seen as they
can be transformed into functional equations of DC type by singling out
the first term of the sum/product, as performed in the examples. 
Not all functions of DC type seem to be elementary, however, for example, 
we were not able to find
an elementary form of $F(z)=\sum 2^{e_0(i)}z^i$, where~$e_0(n)$ counts
the number of zeros in the binary representation of~$n$. On the other hand,
$F(z)$ is of DC type since $F(z^2)=(1+F(z))/(z+2)$.

In the next section, we show that elementary DC functions admit simple 
recurrences, and that many of the well known sequences concerning the binary
representation of~$n$ are generated by them. This will likely allow to
handle whole classes of sequences of such type analytically --- Prodinger
has shown\cite{b-prod} that applying the Mellin transform to power series
generating functions, in order to find asymptotic behaviour, yields useful 
results. Even if this turns out to be difficult in the end, 
the classification of these sequences 
according to their ordinary g.f.s will give identities and clarify the
sequence zoo.

\section{Main Theorem}
Where the theorems in this section are given unproven, there is 
overwhelming evidence for them from our empirical studies --- usually,
generated sequences were compared with recursions up to index 100 and
above. We can now only encourage the reader to find more formal proofs.
We also provide in the next section tables listing the well known sequences
falling into the given categories.
\bigskip\hrule\smallskip
\begin{theo}
Let $\alpha,c,d$ be integers. The following generating functions
of elementary divide-and-conquer type have coefficients satisfying the recurrences

\begin{flalign}
F(z)=\sum_{i=0}^\infty a_iz^i=
\sum_{k=0}^\infty\frac{c^kz^{2^k}}{1-z^{2^k}} & \eq
\left(\begin{tabular}{l}$a_{2n}=ca_n+1$\\$a_{2n+1}=1$\end{tabular}\right),
\quad |c|>0,\label{t1}\\
\sum\frac{c^kz^{2^k}}{1-z^{2^{k+1}}} & \eq
\left(\begin{tabular}{l}$a_{2n}=ca_n$\\$a_{2n+1}=1$\end{tabular}\right),
\quad |c|>0,\label{t2}\\
\prod 1+cz^{2^k} & \eq
\left(\begin{tabular}{l}$a_1=1$\\$a_{2n}=a_n$\\$a_{2n+1}=ca_n$\end{tabular}
\right),
\quad |c|>0,\label{t3}\\
\frac{1}{1-z}\sum\frac{\alpha^k(dz^{2^k}+cz^{2^{k+1}})}{1+z^{2^k}}
& \eq
\left(\begin{tabular}{l}$a_0=0$\\$a_{2n}=\alpha\cdot a_n+c$\\$a_{2n+1}=\alpha\cdot a_{n}+d$\end{tabular}\right),\quad|\alpha|>0.\label{t4}\\
\intertext{Let $c_1,c_2,\ldots,c_D$ be integers, then also}
\prod\left(1+cz^{2^k}+\sum_{i=1}^Dc_iz^{2^{k+1}i}\right) & \eq
\left(\begin{tabular}{l}$a_{n<0}=0,\,a_0=1$\\$a_{2n}=a_n+\sum_{i=1}^Dc_ia_{n-i}$\\$a_{2n+1}=ca_n$\end{tabular}\right),\label{t5}\\
\sum\frac1{1-\sum_{i=1}^Dc_iz^{2^{k}i}} & \eq
\left(\begin{tabular}{l}$a_{2n}=a_n+b_{2n}$\\$a_{2n+1}=b_{2n+1}$\end{tabular}\right),\label{t6}
\end{flalign}
where $b_n$ is a sequence satisfying the linear recurrence 
$b_n=\sum_{i=1}^Dc_ib_{n-i}$.
\end{theo}
\medskip\hrule\bigskip

First, the infinite sums/products are not really; it suffices in practice
to compute the first $\lceil\log_2n\rceil$ terms because any exponent of~$z$
will then be greater than the index of the coefficient we want to get.
This fact can be applied in the proofs, as well. A symbolic calculator
program like PARI helps with visualizing what subsequences are generated.

To prove \eqref{t1}, see that $\frac1{1-z^e}$ generates 
$a_n=[\,e \text{\ divides\ } n\,]$,
and any factor $z^e$ shifts the sequence to the right by $e$ places.
Therefore, 
$$\sum v_2(i)z^i=\sum_{k\ge1}\frac{z^{2^k}}{1-z^{2^k}},\quad
\text{where\ }v_2(n)=-1+\sum\;[\,2^k \text{\ divides\ } n\,],$$
the exponent of the highest power of~$2$ dividing~$n$. Consequently, the
g.f.~in \eqref{t1} adds~$c^k$ if~$2^k$ divides~$n$. We then have 
$$a_n=\sum_{k=0}^{v_2(n)} c^k \quad=\quad
\begin{cases}1+c\sum_{k=0}^{v_2(n)-1}c^k, & \text{for even~$n$;}\\
1,& \text{for odd~$n$.}\end{cases}$$
Similar observations apply to~\eqref{t2}.

Although~\eqref{t3} is a special case of~\eqref{t5}, it deserved to be
mentioned, as there is a simplified form for the coefficients:
$a_n=c^{e_1(n)}$. It is obvious that every exponent~$\le2^k$ is generated
by the first~$k$ terms of the product, and that, every time a one appears
in the binary representation of~$n$, the coefficient of~$z^n$ is multiplied
with~$c$. Finally, the recurrence is an automaton that has the same
behaviour, and the recurrence for~$e_1(n)$ shows it too, just with
addition instead of multiplication.

The principle of a machinery that does one thing when we have a one bit,
and the other with a zero bit, is exploited to the fullest with~\eqref{t4},
where we have a multiplication with~$\alpha$ for any bit, followed by
addition of~$c$ or~$d$, respectively. But this means also that any
recurrence/g.f.~of the form~\eqref{t4} is a linear combination
of two functions, where we have for the case~$\alpha=1$
$$\sum e_1(i)z^i=\frac1{1-z}\sum_{k\ge0}\frac{z^{2^k}}{1+z^{2^k}}, \qquad
\sum e_0(i)z^i=\frac1{1-z}\sum_{k\ge0}\frac{z^{2^{k+1}}}{1+z^{2^k}}.$$
We will list the function pairs more completely in the examples section.
\hfill\qed\break\medskip

As said, \eqref{t5} and~\eqref{t6} are given as well supported
conjectures. It might be possible to construct from them o.g.f.s
of several classes of DC type sequences where only the recurrences
are known.

For another generalization of functions dependent on the binary
representation, the reader is referred to Dumas who 
discusses\cite{b-dumas-www} 2-rational sequences that are defined using linear algebra.
One defines a {\bf 2-rational sequence} $u_n$ by a linear
representation.  Such a  representation consists of a row matrix
$\lambda$, a column matrix $\gamma$, and two square matrices
$A_0$ and $A_1$ whose sizes are respectively $1\times N$, 
$N\times 1$, and $N\times N$ for a certain integer $N$.  If
an integer $n$ has as  binary expansion
\begin{equation*}
n=(n_{\ell}\ldots n_0)_2,
\end{equation*}
the value of the sequence for $n$ is
\begin{equation}
u_n=\lambda A_{n_{\ell}}\dotsb A_{n_0}\gamma.
\end{equation}

\section{Example Sequences}
The use of Neil Sloane's \emph{Online Encyclopedia of Integer Sequences}\cite{b-oeis}
was invaluable in finding the given sequences, and we provide their
A-numbers from the database for further references.

\bigskip\begin{tabular}{llll}\toprule
Type & Parameter(s) & Name/Description & OEIS/References \\\midrule
\eqref{t1} & $c=1$ & $2^{a_n}$ divides $2n, v_2(n)+1$ & A001511, \cite{b-kreg}\\
 & $c=2$ & $n$ XOR $n-1$ & A038712\\
 & $c=-1$ & first Feigenbaum symbolic seq. & A035263, \cite{b-kohel}\\
 \midrule
\eqref{t2} & $c=2$ & highest power of 2 dividing~$n$ & A006519\\
 \midrule
\eqref{t3} & $c=2$ & Gould's seq., $2^{e_1(n)}$ & A001316, \cite{b-kreg}\\
 & $c=3$ & $3^{e_1(n)}$ & (A048883), \cite{b-pt}\\
 & $c=-1$ & Thue-Morse seq.~on $\{1,-1\}$ & \cite{b-astm}\\
\bottomrule\end{tabular}

\begin{tabular}{llll}\toprule
Type & Parameter(s) & Name/Description & OEIS/References \\\midrule
\eqref{t4} & $c=0,d=1$ & ones-counting seq., $e_1(n),\nu(n)$ & A000120, \cite{b-kreg}\\
$\alpha=1$ & $c=1,d=0$ & $e_0(n)$ & A023416\\
 & $c=1,d=1$ & binary length & A070939\\
 & $c=1,d=-1$ & $e_0(n)-e_1(n)$ & A037861\\
 & $c=2,d=1$ & a stopping problem & (A061313)\\
 \midrule
$\alpha=2$ & $c=0,d=1$ & natural numbers, $n$ & A000027\\
 & $c=1,d=0$ & interchange 0s and 1s & A035327, \cite{b-kreg2}\\
 & $c=1,d=1$ & $a_{n-1}$ OR $n$ & A003817 \\
 \midrule
$\alpha=-1$ & $c=0,d=1$ & alternating bit sum for $n$ & A065359, \cite{b-kreg2}\\
 & $c=1,d=0$ & & A083905\\
 & $c=1,d=1$ & & A030300\\
 \midrule
$\alpha=3$ & $c=0,d=1$ & ternary$(n)$ contains no 2 & A005836, \cite{b-kreg}\\
$\alpha=4$ & $c=0,d=1$ & Moser--de Bruijn sequence & A000695, \cite{b-kreg}\\
 \midrule
\eqref{t5} & $c=1,c_1=1$ & Stern-Brocot (Carlitz) seq. & A002487,\cite{b-kreg}\cite{b-kreg2}\\
 & $c=1,c_1=-1$ & a fractal sequence & A005590, \cite{b-kreg}\\
 & $c=3,c_1=2$ & odd entries in Pascal($1\ldots n$) & A006046
\\\bottomrule\end{tabular}

\section{Questions}
What is $\sum2^{e_0(i)}z^i$? And what function generates Per N\o{}rg\aa{}rd's
infinity sequence, defined as $\{a_0=0, a_{2n}=-a_n, a_{2n+1}=a_n+1\}$?
Can the sequences in this work expressed in 2-rational form?

\bigskip
\bigskip


\end{document}